\documentclass[12pt]{amsart}
\usepackage{latexsym,amsmath,amsfonts,amssymb,amsthm}
\theoremstyle{plain}

\newtheorem{The}{Theorem}[section]
\newtheorem{Lem}{Lemma}[section]

\theoremstyle{definition}
\newtheorem*{Ack}{Acknowledgment}
\theoremstyle{remark}

\def\ord{{\rm ord}}
\def\N{\mathbb{N}}
\def\Z{\mathbb{Z}}
\def\Q{\mathbb{Q}}

\def\pmod #1{\ ({\rm{mod}}\ {#1})}
\newcommand{\eu}[2]{\genfrac{\langle}{\rangle}{0pt}{}{#1}{#2}}

\begin{document}
\title{Congruences on Stirling numbers and Eulerian numbers}
\author{Hui-Qin Cao}
\author{Hao Pan}
\address{Department of Mathematics, Nanjing University,
Nanjing 210093, People's Republic of China}
\email{caohq@nau.edu.cn}
\email{haopan79@yahoo.com.cn} 
\date{}
\maketitle

\begin{abstract}
In this paper, we establish some Fleck-Weisman type and Davis-Sun type congruences for the Stirling numbers and the Eulerian numbers.
\end{abstract}

\section{Introduction}
\setcounter{equation}{0}
As usual, we set $\binom{x}{0}=1$ and
$$
\binom{x}{k}=\frac{x(x-1)\cdots(x-k+1)}{k!}\quad\text{for}\quad
k=1,2,\ldots.
$$
We also set $\binom{x}{k}=0$ for any negative integer $k$.

Let $p$ be a prime, and let $n>0, r$ be integers. In 1913, A.
Fleck (cf. \cite{Di}, p. 274) discovered that
\begin{equation}
\sum_{k\equiv r\pmod{p}}\binom{n}{k}(-1)^k\equiv
0\pmod{p^{\lfloor\frac{n-1}{p-1}\rfloor}},
\end{equation}
where $\lfloor\cdot\rfloor$ is the floor function. In 1977, C. S.
Weisman \cite{We} extended Fleck's congruence to prime power
moduli in the following way:
\begin{equation}
\sum_{k\equiv r\pmod{p^{\alpha}}}\binom{n}{k}(-1)^k\equiv
0\pmod{p^{\lfloor\frac{n-p^{\alpha-1}}{p^{\alpha-1}(p-1)}\rfloor}},
\end{equation}
where $\alpha$ is a positive integer.

In 2005, in his lecture notes on Fontaine's rings, D. Wan got another extension of Fleck's
congruence:
\begin{equation}
\label{w1} \sum_{k\equiv
r\pmod{p}}\binom{n}{k}(-1)^k\binom{(k-r)/p}{l}\equiv
0\pmod{p^{\lfloor\frac{n-lp-1}{p-1}\rfloor}},
\end{equation}
where $l\in\N$ and $n>lp$. Later, by a combinatorial approach, Z.
W. Sun \cite{Su1} established a common generalization of Weisman's
and Wan's extensions of Fleck's congruence:
\begin{equation}
\label{su1} \ord_p\bigg(\sum_{k\equiv
r\pmod{p^{\beta}}}\binom{n}{k}(-1)^k\binom{\lfloor(k-r)/p^{\alpha}\rfloor}{l}\bigg)\geq
\left\lfloor\frac{n-p^{\alpha-1}-l}{p^{\alpha-1}(p-1)}\right\rfloor-(l-1)\alpha-\beta
\end{equation}
provided that $\alpha\geq\beta\geq 0$ and $n\geq p^{\alpha-1}$, where
$\ord_p(a)=\sup\{i\in\N : p^i\mid a\}$
is the $p$-adic order of $a\in\Z$.

In fact, with help of $\psi$-operator in Fontaine's theory of
$(\phi,\Gamma)$-modules, (\ref{w1}) and (\ref{su1}) can be improved as
follows \cite{Wa}:
\begin{equation}
\label{w2} \ord_p\bigg(\sum_{k\equiv
r\pmod{p^{\alpha}}}\binom{n}{k}(-1)^k\binom{(k-r)/p^{\alpha}}{l}\bigg)\geq
\left\lfloor\frac{n-p^{\alpha-1}-lp^{\alpha}}{p^{\alpha-1}(p-1)}\right\rfloor.
\end{equation}
And a combinatorial proof of (\ref{w2}) is given in \cite{SW1}. On
the other hand, motivated by algebraic topology, D. M. Davis and
Z. W. Sun \cite{DS, SD} showed that
\begin{equation}
\label{ds} \ord_p\bigg(\sum_{k\equiv
r\pmod{p^{\alpha}}}\binom{n}{k}(-1)^k\binom{(k-r)/p^{\alpha}}{l}\bigg)\geq
\ord_p\bigg(\bigg\lfloor\frac{n}{p^\alpha}\bigg\rfloor!\bigg)
\end{equation}
and
\begin{equation}
\label{sd} \ord_p\bigg(\sum_{k\equiv
r\pmod{p^{\alpha}}}\binom{n}{k}(-1)^k\binom{(k-r)/p^{\alpha}}{l}\bigg)\geq
\ord_p\bigg(\bigg\lfloor\frac{n}{p^{\alpha-1}}\bigg\rfloor!\bigg)-l-\ord_p(l!).
\end{equation}
Note that (\ref{ds}) and (\ref{sd}) can't be deduced from (\ref{w2}). For the
further developments on (\ref{w2}) and (\ref{ds}), the reader may
refer to \cite{Su2, SW1, SW2, Su3}.

The Stirling number $s(n,k)$ of the first kind denotes the number of
permutations of $\{1,2,\ldots,n\}$ which contain exactly $k$
permutation cycles. $s(n,k)\ (0\leq k\leq n)$ can be given by
$$
x(x+1)\cdots(x+n-1)=\sum_{k=0}^{n}s(n,k)x^k.
$$
Similarly, the Stirling number $S(n,k)\ (k\in\N)$ of the second
kind is the number of ways to partition a set of cardinality
$n$ into $k$ nonempty subsets. It is well known that
$$
x^n=\sum_{k=0}^{n}S(n,k)k!\binom{x}{k}
$$
for $n\in\N$. In particular, we set $s(0,0)=S(0,0)=1$ and
$s(n,k)=S(n,k)=0$ whenever $k>n$.

The Eulerian numbers are another special numbers related to
permutations. For an arbitrary permutation $\pi=a_1a_2\cdots a_n$
of $\{1,2,\ldots,n\}$, we say that an element
$i\in\{1,2,\ldots,n-1\}$ is an ascent of $\pi$ if $a_i<a_{i+1}$.
The Eulerian number $\eu nk$ is the number of permutations of
$\{1,2,\ldots,n\}$ having $k$ ascents (cf. \cite{GKP}, p. 267).
(Another commonly used notation is $A(n,k)$ (sometimes $A_{n,k}$)
with $A(n,k)=\eu{n}{k-1}$.) Clearly $\eu n0=1$ and $\eu nk=0$ for
every $k>n-1$. We also set $\eu nk=0$ when $k<0$. It is easy to
check that the Eulerian numbers satisfy the recurrence relation
\begin{align*}
\eu{n}{k}=(k+1)\eu{n-1}{k}+(n-k)\eu{n-1}{k-1}.
\end{align*}

Stirling numbers and Eulerian numbers play important roles in
enumerative combinatorics. In this paper, motivated by (\ref{w2})
and (\ref{sd}), we shall give some similar congruences for the Stirling numbers and the
Eulerian numbers. Firstly, we have the following result for the Eulerian numbers.
\begin{The}
\label{mt1}
Let $p$ be a prime. Let $n>0, r$ be integers. Then for positive
integer $\alpha$ and $l\in\N$, we have
\begin{equation}
\label{ec1}
\ord_p\bigg(\sum_{k\equiv r\pmod{p^{\alpha}}}\eu nk
\binom{(k-r)/p^{\alpha}}{l}\bigg)\geq\ord_p\bigg(\left\lfloor\frac{n}{p^{\alpha-1}}\right\rfloor!\bigg)-\left\lceil
\frac{p^{\alpha-1}+lp^{\alpha}}{p^{\alpha-1}(p-1)}\right\rceil,
\end{equation}
where $\lceil\cdot\rceil$ is the
ceiling function. Moreover, if $a$ is an integer with $a\equiv 1\pmod{p}$,
then
\begin{equation}
\label{ec2}
\sum_{k\equiv r\pmod{p^{\alpha}}}\eu nk a^{k}\equiv
0\pmod{p^{\ord_p(\lfloor n/p^{\alpha-1}\rfloor!)-1}}
\end{equation}
provided that $n\geq p^\alpha$.
\end{The}
The results on Stirling numbers are a little complicated.
\begin{The}
\label{mt2}
Let $p$ be a prime and $n, m$ be positive integers. For arbitrary
integers $a$ and $r$,
\begin{equation}
\label{sc1} \ord_p\bigg(\sum_{k\equiv
r\pmod{p-1}}s(n,k)S(k,m)a^k\bigg)\geq\ord_p(n!)-\ord_p(m!).
\end{equation}
Moreover, if $f(x)$ is a polynomial with integral coefficients,
then
\begin{equation}
\label{sc2} \ord_p\bigg(\sum_{k\equiv
r\pmod{p-1}}s(n,k)f(k)a^k\bigg)\geq \ord_p(n!)-\log_p\binom{n}{l},
\end{equation}
where $l=\min\{\deg f, \lfloor n/p\rfloor\}$
\end{The}
Also, we have the following Weisman type congruence.
\begin{The}
\label{mt3}
Let $p$ be a prime and $n, m$ be positive integers. If $\alpha$ is a positive integer, then
\begin{equation}
\label{sc3} \ord_p\bigg(\sum_{k\equiv
r\pmod{p^\alpha(p-1)}}s(n,k)S(k,m)a^k\bigg)\geq
\bigg\lfloor\frac{n-p^\alpha}{p^\alpha(p-1)}\bigg\rfloor-\ord_p(m!)
\end{equation}
for any integers $a$ and $r$.
\end{The}
The proofs of Theorems \ref{mt1}-\ref{mt3} will be given in the next sections.

\section{Proofs of Theorems \ref{mt1}}
\setcounter{equation}{0}
\begin{Lem}
\label{ecl1}
Let $p$ be a prime and let $n, k\in\N$.
Then for any positive integer $\alpha$,
\begin{equation}
\label{S4}
\ord_p(k!S(n,k))\geq\ord_p\bigg(\left\lfloor\frac{n}{p^{\alpha-1}}\right\rfloor!\bigg)-\left\lfloor\frac{n-k}{p^{\alpha-1}(p-1)}\right\rfloor.
\end{equation}
\end{Lem}
\begin{proof}
We use an induction on $n$. There is nothing to do when $n=0$.
Below we assume that $n\geq 1$ and (\ref{S4}) holds for smaller
values of $n$. Obviously (\ref{S4}) holds for $k=0$. Suppose that
$k\geq 1$. It is known (cf. \cite{Co}, p. 209) that
\begin{equation}
\label{S3}
k!S(n,k)=\sum_{i=k-1}^{n-1}\binom{n}{i}(k-1)!S(i,k-1)\qquad(k\geq1).
\end{equation}
Observe that
\begin{align*}
\ord_p\bigg(\binom{n}{i}\bigg)=&\sum_{j=1}^{\infty}
\bigg(\left\lfloor\frac{n}{p^j}\right\rfloor-\left\lfloor\frac{n-i}{p^j}\right\rfloor-\left\lfloor\frac{i}{p^j}\right\rfloor\bigg)\\
\geq &\sum_{j=\alpha}^{\infty}
\bigg(\left\lfloor\frac{n}{p^j}\right\rfloor-\left\lfloor\frac{n-i}{p^j}\right\rfloor-\left\lfloor\frac{i}{p^j}\right\rfloor\bigg)\\
=&\ord_p\bigg(\left\lfloor\frac{n}{p^{\alpha-1}}\right\rfloor
!\bigg)-\ord_p\bigg(\left\lfloor\frac{n-i}{p^{\alpha-1}}\right\rfloor
!\bigg)-\ord_p\bigg(\left\lfloor\frac{i}{p^{\alpha-1}}\right\rfloor!\bigg).
\end{align*}
By the induction hypothesis, for $k-1\leq i\leq n-1$ we have
\begin{align*}
&\ord_p\bigg(\binom{n}{i}(k-1)!S(i,k-1)\bigg)\\
\geq &\ord_p\bigg(\binom{n}{i}\bigg)+\ord_p\bigg(\left\lfloor\frac{i}{p^{\alpha-1}}\right\rfloor!\bigg)-\left\lfloor\frac{i-(k-1)}{p^{\alpha-1}(p-1)}\right\rfloor\\
\geq&\ord_p\bigg(\left\lfloor\frac{n}{p^{\alpha-1}}\right\rfloor!\bigg)
-\ord_p\bigg(\left\lfloor\frac{n-i}{p^{\alpha-1}}\right\rfloor!\bigg)
   -\left\lfloor\frac{i-k+1}{p^{\alpha-1}(p-1)}\right\rfloor\\
\geq &\ord_p\bigg(\left\lfloor\frac{n}{p^{\alpha-1}}\right\rfloor!\bigg)-\left\lfloor\frac{n-i-p^{\alpha-1}}{p^{\alpha-1}(p-1)}\right\rfloor-\left\lfloor\frac{i-k+1}{p^{\alpha-1}(p-1)}\right\rfloor\\
\geq
&\ord_p\bigg(\left\lfloor\frac{n}{p^{\alpha-1}}\right\rfloor!\bigg)-\left\lfloor\frac{n-k+1-p^{\alpha-1}}{p^{\alpha-1}(p-1)}\right\rfloor\\
\geq
&\ord_p\bigg(\left\lfloor\frac{n}{p^{\alpha-1}}\right\rfloor!\bigg)-\left\lfloor\frac{n-k}{p^{\alpha-1}(p-1)}\right\rfloor.
\end{align*}
This concludes our proof.
\end{proof}
\begin{Lem}
\label{ecl2}
Let $n$ be a positive integer. Then for any polynomial $f(x)$ we have
\begin{equation}
\label{E1}
\sum_{k}\eu nk
f(k)x^k=\sum_{m}m!S(n,m)\sum_{i=0}^{n-m}\binom{n-m}{i}(-1)^{n-m-i}f(i)x^i
\end{equation}
\end{Lem}
\begin{proof}
It is sufficient to prove (\ref{E1}) for $f(x)=x^l$,
$l\in\mathbb{N}$. In the case $l=0$, (\ref{E1}) reduces to
\begin{align}
\label{E2}
 \sum_{k}\eu nk x^k
&=\sum_{m}m!S(n,m)\sum_{i=0}^{n-m}\binom{n-m}{i}(-1)^{n-m-i}x^i\nonumber\\
&=\sum_{m}m!S(n,m)(x-1)^{n-m}.
\end{align}
It is true (cf. \cite{GKP}, p. 269). Now assume that $l>0$
and (\ref{E1}) holds for $l-1$, that is,
\begin{align*}
\sum_{k}\eu nk k^{l-1}x^k
=\sum_{m}m!S(n,m)\sum_{i=0}^{n-m}\binom{n-m}{i}(-1)^{n-m-i}i^{l-1}x^i.
\end{align*}
It follows that
\begin{align*}
\sum_{k}\eu nk k^{l}x^{k-1}
=\sum_{m}m!S(n,m)\sum_{i=1}^{n-m}\binom{n-m}{i}(-1)^{n-m-i}i^{l}x^{i-1}.
\end{align*}
Therefore (\ref{E1}) holds for $l$. We are done.
\end{proof}
\begin{proof}[Proof of (\ref{ec1})]
Let $\zeta$ be a primitive $p^{\alpha}$-th root of the unity. Note
that $\binom{(x-r)/p^{\alpha}}{l}$ is a polynomial in $x$ with the
degree $l$. By Lemma \ref{ecl2}, we have
\begin{align*}
&\sum_{k\equiv r\pmod{p^{\alpha}}}\eu nk \binom{(k-r)/p^{\alpha}}{l}\\
=&\sum_{k}\eu nk \binom{(k-r)/p^{\alpha}}{l}\frac{1}{p^{\alpha}}\sum_{j=0}^{p^{\alpha}-1}\zeta^{j(k-r)}\\
=&\frac{1}{p^{\alpha}}\sum_{j=0}^{p^{\alpha}-1}\zeta^{-jr}\sum_{k}\eu nk \binom{(k-r)/p^{\alpha}}{l}\zeta^{jk}\\
=&\frac{1}{p^{\alpha}}\sum_{j=0}^{p^{\alpha}-1}\zeta^{-jr}\sum_{m}m!S(n,m)\sum_{i=0}^{n-m}\binom{n-m}{i}(-1)^{n-m-i}\binom{(i-r)/p^{\alpha}}{l}\zeta^{ji}\\
=&\sum_{m=0}^{n}m!S(n,m)\sum_{i=0}^{n-m}\binom{n-m}{i}(-1)^{n-m-i}\binom{(i-r)/p^{\alpha}}{l}\frac{1}{p^{\alpha}}\sum_{j=0}^{p^{\alpha}-1}\zeta^{j(i-r)}\\
=&\sum_{m=0}^{n}m!S(n,m)\sum_{i\equiv
r\pmod{p^{\alpha}}}\binom{n-m}{i}(-1)^{n-m-i}\binom{(i-r)/p^{\alpha}}{l}
\end{align*}
Applying Lemma \ref{ecl1} and (\ref{w2}), for every $0\leq m\leq n$ we
have
\begin{align*}
&\ord_p\bigg(m!S(n,m)\sum_{i\equiv
r\pmod{p^{\alpha}}}\binom{n-m}{i}(-1)^{n-m-i}\binom{(i-r)/p^{\alpha}}{l}\bigg)\\
\geq
&\ord_p\bigg(\left\lfloor\frac{n}{p^{\alpha-1}}\right\rfloor!\bigg)-\left\lfloor\frac{n-m}{p^{\alpha-1}(p-1)}\right\rfloor+\left\lfloor\frac{n-m-p^{\alpha-1}-lp^{\alpha}}{p^{\alpha-1}(p-1)}\right\rfloor\\
\geq
&\ord_p\bigg(\left\lfloor\frac{n}{p^{\alpha-1}}\right\rfloor!\bigg)-\left\lceil
\frac{p^{\alpha-1}+lp^{\alpha}}{p^{\alpha-1}(p-1)}\right\rceil.
\end{align*}
This concludes the proof.
\end{proof}
\begin{proof}[Proof of (\ref{ec2})]
Let $\zeta$ be a primitive $p^{\alpha}$-th root of the unity. Using
(\ref{E2}) we have
\begin{align*}
\sum_{k\equiv r\pmod{p^{\alpha}}}\eu nk a^{k}
=&\sum_{k}\eu nk a^{k}\frac{1}{p^{\alpha}}\sum_{j=0}^{p^{\alpha}-1}\zeta^{j(k-r)}\\
=&\frac{1}{p^{\alpha}}\sum_{j=0}^{p^{\alpha}-1}\zeta^{-jr}\sum_{k}\eu nk (a\zeta^{j})^{k}\\
=&\frac{1}{p^{\alpha}}\sum_{j=0}^{p^{\alpha}-1}\zeta^{-jr}\sum_{m}m!S(n,m)(a\zeta^{j}-1)^{n-m}\\
=&\frac{1}{p^{\alpha}}\sum_{j=0}^{p^{\alpha}-1}\zeta^{-jr}\sum_{m=0}^{n}m!S(n,m)\sum_{i=0}^{n-m}\binom{n-m}{i}(-1)^{n-m-i}a^{i}\zeta^{ji}\\
=&\sum_{m=0}^{n}m!S(n,m)(-1)^{n-m}\sum_{i=0}^{n-m}\binom{n-m}{i}(-a)^{i}\frac{1}{p^{\alpha}}\sum_{j=0}^{p^{\alpha}-1}\zeta^{j(i-r)}\\
=&\sum_{m=0}^{n}m!S(n,m)(-1)^{n-m}\sum_{i\equiv
r\pmod{p^{\alpha}}}\binom{n-m}{i}(-a)^{i}.
\end{align*}
In view of Lemma \ref{ecl1} and (1.8) of \cite{Su1},
\begin{align*}
&\ord_p\bigg(m!S(n,m)\sum_{i\equiv
r\pmod{p^{\alpha}}}\binom{n-m}{i}(-a)^{i}\bigg)\\
\geq
&\ord_p\bigg(\left\lfloor\frac{n}{p^{\alpha-1}}\right\rfloor!\bigg)-\left\lfloor\frac{n-m}{p^{\alpha-1}(p-1)}\right\rfloor+\left\lfloor\frac{n-m-p^{\alpha-1}}{p^{\alpha-1}(p-1)}\right\rfloor\\
\geq
&\ord_p\bigg(\left\lfloor\frac{n}{p^{\alpha-1}}\right\rfloor!\bigg)-1
\end{align*}
for every $0\leq m\leq n$. This ends the proof.
\end{proof}

\section{Congruences for Stirling numbers: $k\equiv r\pmod{p-1}$}
\setcounter{equation}{0}
In this section, we shall prove Theorem \ref{mt2}.
Let $\Z_p$ (resp. $\Q_p$) denote the rational $p$-adic integers
ring (resp. field).
\begin{Lem}
\label{scl1}
For any $n\in\N$ and $x\in\Z_p$, $\binom{x}{n}$ is $p$-integral.
\end{Lem}
\begin{proof}
We may choose $x'\in\N$ such that
$$
x\equiv x'\pmod{p^{\ord_p(n!)+1}}.
$$
Then
$$
\binom{x}{n}=\frac{x(x-1)\cdots(x-n+1)}{n!}\equiv \frac{x'(x'-1)\cdots(x'-n+1)}{n!}=\binom{x'}{n}\pmod{p}.
$$
This concludes that $\binom{x}{n}\in\Z_p$ since $\binom{x'}{n}\in\Z$.
\end{proof}
\begin{proof}[Proof of (\ref{sc1})]
Let $\omega$ be the Teichm\"uller character of the multiplicative
group $(\Z/p\Z)^*$. We know that
$$
\omega(a)\in\Z_p,\qquad a=1,2,\ldots,p-1
$$
are exactly all $(p-1)$-th roots of unity in $\Q_p$. And if $g$ is
a primitive root of $p$, then $\omega(g)$ is a $(p-1)$-th
primitive roots of unity. Let $\omega\in\Z_p$ be an arbitrary
$(p-1)$-th primitive roots of unity in $\Q_p$. Thus
\begin{align*}
&\sum_{k\equiv r\pmod{p-1}}s(n,k)S(k,m)a^k\\
=&\sum_{k=0}^ns(n,k)S(k,m)a^k\cdot\frac{1}{p-1}\sum_{j=1}^{p-1}\omega^{j(k-r)}\\
=&\frac{1}{p-1}\sum_{j=1}^{p-1}\omega^{-jr)}\sum_{k}s(n,k)S(k,m)\big(a\omega^j\big)^k\\
=&\frac{1}{p-1}\sum_{j=1}^{p-1}\omega^{-jr}\sum_{k}s(n,k)\bigg(\frac{1}{m!}\sum_{i=0}^{m}\binom{m}{i}(-1)^{m-i}i^k\bigg)\big(a\omega^j\big)^k\\
=&\frac{1}{m!(p-1)}\sum_{j=1}^{p-1}\omega^{-jr}\sum_{i=0}^{m}\binom{m}{i}(-1)^{m-i}\sum_{k}s(n,k)\big(ai\omega^j\big)^k\\
=&\frac{n!}{m!(p-1)}\sum_{j=1}^{p-1}\omega^{-jr}\sum_{i=0}^{m}\binom{m}{i}(-1)^{m-i}\binom{ai\omega^j+n-1}{n}.
\end{align*}
Therefore applying Lemma \ref{scl1},
$$
\ord_p\bigg(\sum_{k\equiv r\pmod{p-1}}s(n,k)S(k,m)a^k\bigg)\geq
\ord_p(n!/m!)=\ord_p(n!)-\ord_p(m!).
$$
\end{proof}
\begin{Lem}
\label{scl2}
Let $n$ and $l$  be nonnegative integers. Then
$$
(n-i)\binom{i}{l-1}\leq\binom{n}{l}
$$
for each integer $0\leq i\leq n$.
\end{Lem}
\begin{proof}
It is easy to check that
$$
(n-i)\binom{i}{l-1}\geq(n-i+1)\binom{i-1}{l-1}\Longleftrightarrow i\leq\frac{(l-1)(n+1)}{l}.
$$
Hence when $n\geq l\geq 1$
\begin{align*}
(n-i)\binom{i}{l-1}\leq&(n-\lfloor(l-1)(n+1)/l\rfloor)\binom{\lfloor(l-1)(n+1)/l\rfloor}{l-1}\\
=&\left\lceil\frac{n-l+1}{l}\right\rceil\binom{n-\lceil(n-l+1)/l\rceil}{l-1}\\
\leq&\frac{n}{l}\binom{n-1}{l-1}=\binom{n}{l}.
\end{align*}
\end{proof}
\begin{proof}[Proof of (\ref{sc2})]
We use an induction on $\deg f$. The case $\deg f=0$ follows from
(\ref{sc1}) by setting $m=1$. Below we assume that $\deg f>0$ and (\ref{sc2}) holds for
the smaller values of $\deg f$. It is known (cf. \cite{Co}, p. 215) that
\begin{equation}
\label{ss3}
ks(n,k)=\sum_{i=k-1}^{n-1}\binom{n}{i}(n-i-1)!s(i,k-1)\qquad(k\geq1).
\end{equation}
Write $f(x)=xf_1(x)+c$ with $\deg f_1=\deg f-1$. Then
\begin{align*}
&\sum_{k\equiv r\pmod{p-1}}s(n,k)f(k)a^k\\
=&\sum_{k\equiv r\pmod{p-1}}f_1(k)a^k\sum_{i=k-1}^{n-1}\binom{n}{i}(n-i-1)!s(i,k-1)+c\sum_{k\equiv r\pmod{p-1}}s(n,k)a^k\\
=&\sum_{i=0}^{n-1}\frac{n!}{i!(n-i)}\sum_{k\equiv
r\pmod{p-1}}s(i,k-1)f_1(k)a^k+c\sum_{k\equiv r\pmod{p-1}}s(n,k)a^k\\
=&\sum_{i=0}^{n-1}\frac{an!}{i!(n-i)}\sum_{k\equiv
r-1\pmod{p-1}}s(i,k)f_1(k+1)a^k+c\sum_{k\equiv r\pmod{p-1}}s(n,k)a^k.
\end{align*}
When $i=0$,
\begin{align*}
&\ord_p\bigg(\frac{an!}{0!(n-0)}\sum_{k\equiv
r-1\pmod{p-1}}s(0,k)f_1(k+1)a^k\bigg)\\
\geq&\ord_p(n!)-\ord_p(n)\geq\begin{cases}0=\ord_p(n!)-\log_p\binom{n}{0}&\text{if }n<p,\\
\ord_p(n!)-\log_p
n\geq\ord_p(n!)-\log_p\binom{n}{l}&\text{otherwise}.
\end{cases}
\end{align*}
For every $0<i\leq n-1$, by the induction hypothesis,
\begin{align*}
&\ord_p\bigg(\frac{an!}{i!(n-i)}\sum_{k\equiv
r-1\pmod{p-1}}s(i,k)f_1(k+1)a^k\bigg)\\
\geq&\ord_p(n!)-\ord_p(i!)-\ord_p(n-i)+\ord_p(i!)-\log_p\binom{i}{l'},
\end{align*}
where $l'=\min\{\deg f-1,\lfloor i/p\rfloor\}$. It suffices to show that
$$
\ord_p(n-i)+\log_p\binom{i}{l'}\leq\log_p\binom{n}{l}.
$$
When $i>n-p$, clearly $l-1\leq l'\leq l$ and $\ord_p(n-i)=0$. Hence
$$
\binom{i}{l'}\leq\max\{\binom{n-1}{l-1},\binom{n-1}{l}\}\leq\binom{n}{l}.
$$
Below assume that $i\leq n-p$.
If $\deg f-1\leq i/p$, then applying Lemma \ref{scl2},
$$
(n-i)\binom{i}{\deg f-1}\leq\binom{n}{\deg f}\leq\binom{n}{l}
$$
since $\deg f\leq\lfloor i/p\rfloor+1\leq\lfloor n/p\rfloor$ now.
Also, when $i/p<\deg f-1$, we have
$$
(n-i)\binom{i}{\lfloor i/p\rfloor}\leq\binom{n}{\lfloor i/p\rfloor+1}\leq\binom{n}{l}
$$
provided that $\lfloor i/p\rfloor<\lfloor n/p\rfloor$. In two above cases, we both obtain that
$$
\ord_p(n-i)+\log_p\binom{i}{l'}\leq\log_p(n-i)+\log_p\binom{i}{l'}\leq\log_p\binom{n}{l}.
$$
This concludes our proof.
\end{proof}

\section{Congruences for Stirling numbers: $k\equiv r\pmod{p^\alpha(p-1)}$}
\setcounter{equation}{0}

In this section, we shall prove Theorem \ref{mt3}.
Define
$$
C_{d,r}(n,m,a)=\sum_{k\equiv r\pmod{d}}s(n,k)S(k,m)a^k.
$$
Let $\zeta_d$ be a primitive $d$-th root of the unity. Then
\begin{align}
C_{d,r}(n,m,a)=&\sum_{k}s(n,k)S(k,m)a^k\bigg(\frac{1}{d}\sum_{j=0}^{d-1}\zeta_d^{j(k-r)}\bigg)\notag\\
=&\frac{1}{d}\sum_{j=0}^{d-1}\zeta_d^{-jr}\sum_{k}s(n,k)a^k\frac{1}{m!}\sum_{i=0}^{m}\binom{m}{i}(-1)^{m-i}i^k\zeta_d^{jk}\notag\\
=&\frac{1}{m!d}\sum_{j=0}^{d-1}\zeta_d^{-jr}\sum_{i=0}^{m}\binom{m}{i}(-1)^{m-i}\sum_{k}s(n,k)(ai\zeta_d^{j})^{k}\notag\\
=&\frac{(-1)^n}{m!d}\sum_{i=0}^{m}\binom{m}{i}(-1)^{m-i}\sum_{j=0}^{d-1}\zeta_d^{-jr}(-ai\zeta_d^{j})_{n}.
\end{align}
\begin{Lem}
\label{scl3} Let $p$ be a prime and $\alpha$ be a positive
integer. Then for any $1\leq k\leq p^\alpha(p-1)$, we have
\begin{equation}
\label{scl3e}
s(p^\alpha(p-1),k)\equiv\begin{cases}
1\pmod{p}\quad&\text{if }k\equiv 0\pmod{p^{\alpha-1}(p-1)},\\
0\pmod{p}\quad&\text{otherwise}.
\end{cases}
\end{equation}
\end{Lem}
\begin{proof}
Let $x$ be a. Apparently
$$
x(x+1)\cdots(x+p-1)\equiv x^p-x\pmod{p}.
$$
Thus
\begin{align*}
x(x+1)\cdots(x+p^\alpha(p-1)-1)\equiv&(x(x+1)\cdots(x+p-1))^{p^{\alpha-1}(p-1)}\\
\equiv&(x^p-x)^{p^{\alpha-1}(p-1)}\\
=&\sum_{j=0}^{p^{\alpha-1}(p-1)}\binom{p^{\alpha-1}(p-1)}{j}(-1)^jx^{p^{\alpha-1}(p-1)+(p-1)j}\pmod{p}.
\end{align*}
By the Lucas congruence, we know that
$$
\binom{p^{\alpha-1}(p-1)}{j}\equiv\begin{cases}
\binom{p-1}{j/p^{\alpha-1}}\pmod{p}&\text{if }p^{\alpha-1}\mid j,\\
0\pmod{p}&\text{otherwise}.
\end{cases}
$$
Hence
\begin{align*}
\sum_{k=0}^{p^{\alpha-1}(p-1)}s(p^{\alpha-1}(p-1),k)x^k\equiv&
\sum_{j=0}^{p-1}\binom{p^{\alpha-1}(p-1)}{p^{\alpha-1}j}x^{p^{\alpha-1}(p-1)+p^{\alpha-1}(p-1)j}\\
&\sum_{j=1}^{p}\binom{p^{\alpha-1}(p-1)}{p^{\alpha-1}j}x^{p^{\alpha-1}(p-1)j}\pmod{p},
\end{align*}
which is obviously equivalent to (\ref{scl3e}).
\end{proof}
\begin{proof}[Proof of Theorem \ref{mt3}]
Apparently
\begin{align*}
&C_{d,r}(n+1,m,a)\\
=&\sum_{k\equiv r\pmod{d}}s(n+1,k)S(k,m)a^k\\
=&n\sum_{k\equiv r\pmod{d}}s(n,k)S(k,m)a^k+\sum_{k\equiv r\pmod{d}}s(n,k-1)S(k,m)a^k\\
=&n\sum_{k\equiv r\pmod{d}}s(n,k)S(k,m)a^k+\sum_{k\equiv r-1\pmod{d}}s(n,k)(mS(k,m)+S(k,m-1))a^{k+1}\\
=&nC_{d,r}(n,m,a)+amC_{d,r-1}(n,m,a)+aC_{d,r-1}(n,m-1,a).
\end{align*}
Also observe that
$$
\bigg\lfloor\frac{n-p^\alpha}{p^\alpha(p-1)}\bigg\rfloor
=\bigg\lfloor\frac{\lfloor n/p^\alpha\rfloor-1}{p-1}\bigg\rfloor.
$$
Hence without loss of generality, assume that $p^\alpha$ divides
$n$. We make an induction on $n$. Clearly the cases
$n<p^{\alpha+1}$ is trivial. Let $\zeta$ be a primitive $p^\alpha(p-1)$-th root of unity. Then
\begin{align*}
&(-1)^nm!p^\alpha(p-1)C_{p^\alpha(p-1),r}(n+p^\alpha(p-1),m,a)\\
=&\sum_{i=0}^{m}\binom{m}{i}(-1)^{m-i}\sum_{j=1}^{p^\alpha(p-1)}\zeta^{-jr}(-ai\zeta^{j})_{n+p^\alpha(p-1)}\\
=&\sum_{i=0}^{m}\binom{m}{i}(-1)^{m-i}\sum_{j=1}^{p^\alpha(p-1)}\zeta^{-jr}(-ai\zeta^{j})_{n}(-ai\zeta^{j}-n)_{p^\alpha(p-1)}\\
=&\sum_{i=0}^{m}\binom{m}{i}(-1)^{m-i}\sum_{j=1}^{p^\alpha(p-1)}\zeta^{-jr}(-ai\zeta^{j})_{n}\sum_{k=0}^{p^\alpha(p-1)}s(p^\alpha(p-1),k)(ai\zeta^{j}+n)^k\\
=&\sum_{i=0}^{m}\binom{m}{i}(-1)^{m-i}\sum_{j=1}^{p^\alpha(p-1)}\zeta^{-jr}(-ai\zeta^{j})_{n}\sum_{k=0}^{p^\alpha(p-1)}s(p^\alpha(p-1),k)\sum_{l=0}^k\binom{k}{l}(ai\zeta^{j})^ln^{k-l}\\
=&\sum_{i=0}^{m}\binom{m}{i}(-1)^{m-i}\sum_{k=0}^{p^\alpha(p-1)}s(p^\alpha(p-1),k)\sum_{l=0}^k\binom{k}{l}(ai)^ln^{k-l}\sum_{j=1}^{p^\alpha(p-1)}\zeta^{j(l-r)}(-ai\zeta^{j})_{n}\\
\end{align*}
Applying Lemma \ref{scl3}, we have
\begin{align*}
&m!C_{p^\alpha(p-1),r}(n+p^\alpha(p-1),m,a)\\
=&\sum_{i=0}^{m}\binom{m}{i}(-1)^{m-i}\sum_{k=0}^{p^\alpha(p-1)}s(p^\alpha(p-1),k)\sum_{l=0}^k\binom{k}{l}(ai)^ln^{k-l}
C_{p^\alpha(p-1),r-l}(n,1,ai)\\
\equiv&\sum_{i=0}^{m}\binom{m}{i}(-1)^{m-i}\sum_{k=1}^{p}(ai)^{kp^{\alpha-1}(p-1)}C_{p^\alpha(p-1),r-kp^{\alpha-1}(p-1)}(n,1,ai)\pmod{p^{\lfloor\frac{n-p^\alpha}{p^\alpha(p-1)}\rfloor+1}},
\end{align*}
since $p^\alpha\mid n$ and
$$
C_{p^\alpha(p-1),r-kp^{\alpha-1}(p-1)}(n,1,ai)\equiv 0\pmod{p^{\lfloor\frac{n-p^\alpha}{p^\alpha(p-1)}\rfloor}}
$$
by the induction hypothesis on $n$. When $p\mid a$, clearly now
$$
m!C_{p^\alpha(p-1),r}(n+p^\alpha(p-1),m,a)\equiv0\pmod{p^{\lfloor\frac{n-p^\alpha}{p^\alpha(p-1)}\rfloor+1}}.
$$
And if $p\nmid a$,
\begin{align*}
&m!C_{p^\alpha(p-1),r}(n+p^\alpha(p-1),m,a)\\
\equiv&\sum_{\substack{0\leq i\leq m\\ p\nmid
i}}\binom{m}{i}(-1)^{m-i}\sum_{k=1}^{p}\big((ai)^{kp^{\alpha-1}}\big)^{p-1}C_{p^\alpha(p-1),r-kp^{\alpha-1}(p-1)}(n,1,ai)\\
\equiv&\sum_{\substack{0\leq i\leq m\\ p\nmid
i}}\binom{m}{i}(-1)^{m-i}\sum_{k=1}^{p}C_{p^\alpha(p-1),r-kp^{\alpha-1}(p-1)}(n,1,ai)\\
=&\sum_{\substack{0\leq i\leq m\\ p\nmid
i}}\binom{m}{i}(-1)^{m-i}\sum_{k=1}^{p}\sum_{\substack{l\equiv r-kp^{\alpha-1}(p-1)\\\pmod{p^{\alpha}(p-1)}}}s(n,l)(ai)^k\\
=&\sum_{\substack{0\leq i\leq m\\ p\nmid
i}}\binom{m}{i}(-1)^{m-i}\sum_{l\equiv r\pmod{p^{\alpha-1}(p-1)}}s(n,l)(ai)^k\\
=&\sum_{\substack{0\leq i\leq m\\ p\nmid
i}}\binom{m}{i}(-1)^{m-i}C_{p^{\alpha-1}(p-1),r}(n,1,ai)\pmod{p^{\lfloor\frac{n-p^\alpha}{p^\alpha(p-1)}\rfloor+1}}.
\end{align*}
Thus it suffices to show that
$$
\ord_p C_{p^{\alpha-1}(p-1),r}(n,1,ai)\geq\bigg\lfloor\frac{n-p^\alpha}{p^\alpha(p-1)}\bigg\rfloor+1.
$$
Note that $n\geq p^\alpha$ now. If $\alpha=1$, then by (\ref{sc1})
$$
\ord_p
C_{p-1,r}(n,1,ai)\geq\bigg\lfloor\frac{n}{p}\bigg\rfloor\geq\ord_p(n!)\geq\bigg\lfloor\frac{n+p(p-2)}{p(p-1)}\bigg\rfloor=\bigg\lfloor\frac{n-p}{p(p-1)}\bigg\rfloor+1.
$$
Also if $\alpha\geq 2$, by an induction on $\alpha$, we have
$$
\ord_p C_{p-1,r}(n,1,ai)\geq\bigg\lfloor\frac{\lfloor
n/p^{\alpha-1}\rfloor-1}{p-1}\bigg\rfloor\geq\bigg\lfloor\frac{\lfloor
n/p^\alpha\rfloor-1}{p-1}\bigg\rfloor+1.
$$
All are done.
\end{proof}
\begin{Ack}
We thank our advisor, Professor Zhi-Wei Sun, for his helpful
suggestions on this paper. \end{Ack}

\end{document}